\documentclass{amsart}
\usepackage{amsmath,amsfonts,amssymb,amsthm,enumerate,latexsym,mathrsfs}
\usepackage{graphicx}
\newcounter{nthm}

\newcounter{nthmm}

\newtheorem{thm}{Theorem}

\newtheorem{thm*}[nthm]{Theorem}
\newtheorem{thm**}[nthmm]{Th\'eor\`eme}
\newtheorem{defn}{Definition}
\newtheorem*{defn*}{Definition}
\newtheorem{prop}{Proposition}
\newtheorem{prop*}[nthm]{Proposition}

\newtheorem{conj}{Conjecture}

\newcommand{\N}{\mathbb{N}}
\newcommand{\Z}{\mathbb{Z}}

\newcommand{\Q}{\mathbb{Q}}

\newcommand{\K}{\mathcal{K}}

\newcommand{\funct}[2]{#1 \longrightarrow #2}
\newcommand{\arrows}[3]{\longrightarrow {#1}^{#2}_{#3}}
\newcommand{\restrict}[2]{#1 \upharpoonright #2}
\newcommand{\la}[1]{\stackrel{#1}{\longleftarrow}}

\newcommand{\m}[1]{\textbf{#1}}

\newcommand{\mc}[1]{\widetilde{\textbf{#1}}}

\newcommand{\LO}{\mathrm{LO}}

\newcommand{\Aut}{\mathrm{Aut}}

\newcommand{\Age}{\mathrm{Age}}

\author{Lionel Nguyen Van Th\'e}

\address{Laboratoire d'Analyse, Topologie et Probabilit\'es, Universit\'e d'Aix-Marseille, Centre de Math\'ematiques et Informatique (CMI),
Technop\^{o}le Ch\^{a}teau-Gombert, 39, rue F. Joliot Curie, 13453 Marseille Cedex 13, France}

\email{lionel.nguyen-van-the\@latp.univ-mrs.fr}


\date{December 2014}

\title[Structural Ramsey theory with the KPT correspondence in mind]{A survey on structural Ramsey theory and topological dynamics with the Kechris-Pestov-Todorcevic correspondence in mind}

\begin{document}
\maketitle

\section{Introduction} 

The article \cite{KPT}, published in 2005 by Kechris, Pestov and Todorcevic, established a surprising correspondence between structural Ramsey theory and topological dynamics. As an immediate consequence, it triggered a new interest for structural Ramsey theory. The purpose of the present paper is to present a self-contained survey of the corresponding developments. 

The paper is organized as follows: Section \ref{section:RP} covers the fundamentals of structural Ramsey theory and of Fra\"iss\'e theory. Section \ref{section:flows} introduces the necessary notions from dynamics. Section \ref{section:KPT} presents the Kechris-Pestov-Todorcevic (in short, KPT) correspondence. Section \ref{section:applications} indicates a list of situations where the correspondence was successfully applied. Finally, Section \ref{section:open} contains a list of open problems. Most of the material presented here is based on \cite{KPT} as well as on the habilitation memoir \cite{NVT5}. 

\

\emph{Caution:} The published version of the paper, written in 2013 and submitted in final version in 2014, does not contain updates regarding references or results after July 2013. This current version contains more complete references, but is still far from being exhaustive.

\section{Structural Ramsey theory}

\label{section:RP}

The foundational result of Ramsey theory appeared in 1930. It was proved by Ramsey and can be stated as follows (For a set $X$ and a positive integer $l$, $[X]^l$ denotes the set all of subsets of $X$ with $l$ elements):

\begin{thm}[Ramsey \cite{Ra}]
\label{thm:Ramsey}
For every $l, m \in \N$, there exists $p \in \N$ such that for every set $X$ with $p$ elements, if $[X]^l$ is partitioned into two classes $[X]^l = R \cup B$, then there exists $Y \subset X$ with $m$ elements such that $[Y]^l \subset R$ ou $[Y]^l \subset B$.    

\end{thm}

However, it is only at the beginning of the seventies that the essential ideas behind this theorem crystalized and expanded to structural Ramsey  theory. The goal was then to obtain results similar to Ramsey's theorem in a setting where more structure appears. For example, if $\m H$ is a finite graph, there exists a finite graph $\m K$ with the following property: for every coloring of the edges of $\m K$ in two colors, there exists a finite induced subgraph of $\m K$ isomorphic to $\m H$ where all edges receive the same color. Many other results of the same kind now exist for a wide variety of finite structures, and those are at the center of \cite{KPT}.

\subsection{Fundamentals of Fra\"iss\'e theory.}

\label{subsection:fundamentalsFraisse}


\subsubsection{First order structures}

Let $L=\{R_i :~i \in I\} \cup \{ f_{j}:j\in J\}$ be a fixed language, that is to say a list of symbols to be interpreted later as relations and functions, each symbol having a corresponding integer called its \emph{arity}. The arity of the relation symbol $R_{i}$ is a positive integer $\alpha(i)$ and the arity of each function symbol $f_{j}$ is a non-negative integer $\beta(j)$. Let $\m{A}$ and $\m{B}$ be two $L$-structures (that is, non empty sets $A$, $B$ equipped with relations $R_i ^{\m{A}}\subset A^{\alpha(i)}$ and $R_i ^{\m{B}}\subset B^{\alpha(i)}$ for each $i \in I$ and functions $f_j ^{\m{A}} : \funct{A^{\beta(i)}}{A}$ and $f_j ^{\m{B}} : \funct{B^{\beta(i)}}{B}$ for each $j\in J$). An \emph{embedding}\index{embedding} from $\m{A}$ to $\m{B}$
is an injective map $\pi : \funct{A}{B}$ such that for every $i \in I$, $x_1, \ldots, x_{\alpha(i)} \in A$: \[(x_1, \ldots, x_{\alpha}) \in R_i^{\m{A}} \quad \textrm{iff} \quad (\pi (x_1), \ldots, \pi (x_{\alpha(i)})) \in R_i^{\m{B}},\] and every $j \in J$, $x_1, \ldots, x_{\beta (j)} \in A$: \[\pi(f_{j} ^{\m A}(x_1, \ldots, x_{\beta(j)}) = f_{j} ^{\m B}(\pi (x_1), \ldots, \pi (x_{\beta(j)})).\]

An \emph{isomorphism}\index{isomorphism} from $\m{A}$ to $\m{B}$ is a surjective embedding while an \emph{automorphism}\index{automorphism} of $\m{A}$ is an isomorphism from $\m{A}$ onto itself. Of course,
$\m{A}$ and $\m{B}$ are \emph{isomorphic}\index{isomorphic} when there is an isomorphism from
$\m{A}$ to $\m{B}$. This is written $\m{A} \cong \m{B}$. Finally, the set $\binom{\m{B}}{\m{A}}$ is defined as:
\[
\binom{\m{B}}{\m{A}} = \{ \mc{A} \subset \m{B} : \mc{A} \cong \m{A} \}.
\]

This is the set of \emph{copies} of $\m A$ in $\m B$. Above, the notation $\mc{A} \subset \m{B}$ is used to mean that $\mc A$ is a substructure of $\m B$, i.e. that the underlying set of $\m A$ is contained in the underlying set of $\m B$, and that all relations and functions on $\mc A$ are induced by those of $\m B$. Note however that a subset of $B$ may not support a substructure of $\m B$, but that it always \emph{generates} a substructure of $\m B$ in an obvious way. 

 \subsubsection{Fra\"iss\'e theory}

A structure $\m{F}$ is \emph{ultrahomogeneous}\index{ultrahomogeneous!ultrahomogeneous structure} when every isomorphism between finite substructures of $\m{F}$ can be extended to an automorphism of $\m{F}$. When in addition $\m F$ is countable and every finite subset of $F$ generates a finite substructure of $\m F$ (we say in that case that $\m F$ is \emph{locally finite}), it is a \emph{Fra\"iss\'e structure}.

Let $\m{F}$ be an $L$-structure. The \emph{age} of $\m{F}$\index{age}, denoted $\mathrm{Age}(\m{F})$\index{$\mathrm{Age}(\m{F})$}, is the collection of all finitely generated $L$-structures that can be embedded into $\m{F}$. Observe also that if $\m{F}$ is countable, then $\mathrm{Age}(\m{F})$ contains only countably many isomorphism types. Abusing language, we will say that $\mathrm{Age}(\m{F})$ is countable. Similarly, a class $\mathcal{K}$ of $L$-structures will be said to be countable if it contains only countably many isomorphism types. 

If $\m{F}$ is a Fra\"iss\'e $L$-structure, then observe that $\mathrm{Age}(\m{F})$: 

\begin{enumerate}

\item is countable, 

\item is \emph{hereditary}\index{hereditary}: for every $L$-structure $\m{A}$ and every $\m{B} \in \mathrm{Age}(\m{F})$, if $\m{A}$ embeds in $\m{B}$, then $\m{A} \in \mathrm{Age}(\m{F})$.

\item satisfies the \emph{joint embedding property}\index{joint embedding property}: for every $\m{A}, \m{B} \in \mathrm{Age}(\m{F})$, there is $\m{C} \in \mathrm{Age}(\m{F})$ such that $\m{A}$ and $\m{B}$ embed in $\m{C}$.

\item satisfies the \emph{amalgamation property}\index{amalgamation!amalgamation property} (or is an \emph{amalgamation class}\index{amalgamation!amalgamation class}): for every $\m{A}$, $\m{B} _0$, $\m{B} _1 \in \mathrm{Age}(\m{F})$ and embeddings $f_0 : \funct{\m{A}}{\m{B} _0}$ and $f_1 : \funct{\m{A}}{\m{B}_{1}}$, there is $\m{C} \in \mathrm{Age}(\m{F})$ and embeddings $g_0 : \funct{\m{B} _0}{\m{C}}$, $g_1 : \funct{\m{B} _1}{\m{C}}$ such that $g_0 \circ f_0 = g_1 \circ f_1$. 

\item contains structures of arbitrarily high finite size. 

\end{enumerate}

Any class of finitely generated structures satisfying those five items is called a \emph{Fra\"iss\'e class}. The following theorem, due to Fra\"iss\'e, establishes that every Fra\"iss\'e class is actually the age of a Fra\"iss\'e structure. 

\begin{thm}[Fra\"iss\'e \cite{Fr0}]

\label{thm:Fraisse}
\index{Fra\"iss\'e!theorem}

Let $L$ be a relational signature and let $\mathcal{K}$ be a Fra\"iss\'e class of $L$-structures. Then there is, up to isomorphism, a unique Fra\"iss\'e $L$-structure $\m{F}$ such that $\mathrm{Age}(\m{F}) = \mathcal{K}$. The structure $\m{F}$ is called the \emph{Fra\"iss\'e limit}\index{Fra\"iss\'e!Fra\"iss\'e limit} of $\mathcal{K}$ and denoted $\mathrm{Flim}(\mathcal{K})$\index{$\mathrm{Flim}(\mathcal{K})$}. 

\end{thm}  


\subsubsection{Examples of Fra\"iss\'e classes and Fra\"iss\'e limits}

\begin{enumerate}

\item Linear orders: consider the class of all finite linear orders $\mathcal{LO}$. The language consists of one relational symbol $<$, which is binary (has arity $2$). An element of $\mathcal{LO}$ is of the form $\m A = (A, <^{\m A})$, and is made of a set together with a linear order. The class $\mathcal{LO}$ is a Fra\"iss\'e class, and its Fra\"iss\'e limit is nothing else than the usual linear order $(\Q, <^{\Q})$. 

\item Vector spaces: fix a finite field $F$ and consider the class $\mathcal V _{F}$ of all finite vector spaces over $F$. The relevant language consists of one binary function symbol $+$ and finitely many unary function symbols $M_{\lambda}$ ($\lambda \in F$). In a structure $\m A$, $+$ is interpreted as a group operation on $A$, $M_{\lambda}$ as the scalar multiplication by $\lambda$ for each $\lambda \in F$, and all the usual axioms of vector spaces are satisfied. The class $\mathcal V_{F}$ is a Fra\"iss\'e class, and its limit is the vector space $V_{F}$ of countable dimension over $F$.


\item Graphs: in the undirected case (which is the case we will refer to when we mention graphs without any further indication), the language is made of one binary relation symbol $E$. In a structure, $E$ is interpreted as an irreflexive, symmetric relation. There are several Fra\"iss\'e classes of such objects, but all of them have been classified by Lachlan and Woodrow in \cite{LW}. An example of such a class is the class $\mathcal G$ of all finite graphs. The Fra\"iss\'e limit of $\mathcal G$ is the so-called countable random graph. For directed graphs, the language is made of one binary relation symbol $\leftarrow$ which is interpreted as an irreflexive, antisymmetric binary relation. Fra\"iss\'e classes of finite directed graphs have also been classified, but only much later than graphs. This classification is due to Cherlin in \cite{Ch}. 

\end{enumerate}


\subsubsection{Non-Archimedean Polish groups}

Another remarkable feature of Fra\"iss\'e structures is provided by their automorphism groups. Let $\m F$ be a Fra\"iss\'e structure. Because its underlying set is countable, we may assume that this set is actually $\N$ and the group $\Aut(\m F)$ may be thought of as a subgroup of the permutation group of $\N$. Moreover, if $g$ is a permutation of $\N$ failing to be an automorphism of $\m F$, then there is a finite subset of $\N$ on which this failure is witnessed. Therefore, $\Aut(\m F)$ is a closed subgroup of $S_{\infty}$, the permutation group of $\N$ equipped with the pointwise convergence topology. It turns out that \emph{every} closed subgroup of $S_{\infty}$ arises that way. The class of all closed subgroups of $S_{\infty}$ can also be defined abstractly in several ways: it coincides with the class of all Polish groups that admit a basis at the identity consisting of open subgroups, but also with the class of all Polish groups that admit a compatible left-invariant ultrametric \cite{BK}. Recently, it has been referred to as the class of \emph{non-Archimedean Polish groups} (see \cite{Ke}). It includes all countable discrete groups as well as all profinite groups, but in the sequel, we will mostly concentrate on non locally-compact groups.


%
%

\subsection{The Ramsey property}

\label{subsection:RP}

Throughout this section, $L$ is a fixed language. Let $k\in \N$, and $\m A, \m B, \m C$ be $L$-structures. Recall that the set of all copies of $\m A$ in $\m B$ is the set $$ \binom{\m B}{\m A} = \{ \mc A \subset \m B : \mc A \cong \m A\}.$$ 

The standard arrow partition symbol $$ \m C \arrows{(\m B)}{\m A}{k}$$ is used to mean that for every map $c: \funct{\binom{\m C}{\m A}}{[k]:=\{0,1,...,k-1\}}$, thought as a $k$-coloring of the copies of $\m A$ in $\m C$, there is $\mc B \in \binom{\m C}{\m B}$ such that $c$ is constant on $\binom{\mc B}{\m A}$. 

\begin{defn}

A class $\mathcal{K}$ of $L$-structures has the \emph{Ramsey property}, or is a \emph{Ramsey class}, when $$ \forall k \in \N \quad \forall \m A, \m B \in \mathcal{K} \quad \exists \m C \in \mathcal{K} \quad \m C \arrows{(\m B)}{\m A}{k}.$$

\end{defn}

When $\mathcal{K} = \Age (\m F)$, where $\m F$ is a Fra\"iss\'e structure, this is equivalent, via a compactness argument, 
to: $$ \forall k \in \N \quad \forall \m A, \m B \in \mathcal{K} \quad \m F \arrows{(\m B)}{\m A}{k}.$$

In other words, every finite coloring of the copies of $\m A$ in $\m F$ must be constant on arbitrarily large finite sets. The first example of a Ramsey class is provided by Ramsey's theorem, which states that the class of all finite sets (i.e. structures in the empty language), or equivalently of all finite linear orders, forms a Ramsey class. As indicated previously, the search for Ramsey classes generated a considerable activity in the seventies and in the early eighties. The most significant examples of Ramsey classes which appeared during that period are provided by finite Boolean algebras (Graham-Rothschild, \cite{GR}) and by finite vector spaces over a fixed finite field (Graham-Leeb-Rothschild, \cite{GLR1, GLR2}). However, being Ramsey turns out to be very restrictive, and many natural classes of finite structures do not have the Ramsey property, for example finite equivalence relations, finite graphs, finite relational structures in a fixed language, finite $K_{n}$-free graphs, finite posets, etc...Nevertheless, it appears that those classes are in fact not so far from being Ramsey. In particular, they can be expanded into Ramsey classes simply by adding linear orderings. More details will be given in Section \ref{section:applications}.

%
%

The importance of linear orderings, and more generally of rigidity, in relation to the Ramsey property was realized pretty early. A structure is \emph{rigid} when it admits no non-trivial automorphism. Essentially, all Ramsey classes must be made of rigid structures. Sets, Boolean algebras and vector spaces do not fall into that category, but those being Ramsey is equivalent to some closely related classes of rigid structures being Ramsey. To use the common jargon, rigidity prevents the appearance of Sierpi\'nski type colorings, which do not stabilize on large sets. Let us illustrate this on (simple, loopless) directed graphs: let $\m A$ be a directed edge, and $\m B$ be a directed $3$-cycle. Then no directed graph $\m C$ satisfies $\m C \arrows{(\m B)}{\m A}{2}$. Indeed, take a linear ordering $<$ on $C$. Given an edge $x\la{} y$ in $\m C$, color it blue if $x<y$ and red otherwise. Then every copy of $\m B$ in $\m C$ has edges of each color. 

Another restriction imposed by the Ramsey property appears in the following result.  

\begin{prop}[Ne\v{s}et\v{r}il-R\"odl \cite{NR1}, p.294, Lemma 1]

Let $\K$ be a class of finite $L$-structures consisting of rigid elements. Assume that $\K$ has the hereditarity property, the joint embedding property, and the Ramsey property. Then $\K$ has the amalgamation property. 
\end{prop}

This result explains why structural Ramsey theory and Fra\"iss\'e theory are so closely related: when a class of finite structures satisfies very common properties, it has to be Fra\"iss\'e whenever it is Ramsey. Amalgamation itself is a very restrictive feature, and was at the center of a very active area of research in the eighties. In particular, it led to spectacular classification results, the most significant ones being probably those we already mentioned concerning finite graphs (Lachlan-Woodrow, \cite{LW}), finite tournaments (Lachlan \cite{La}, based on the work of Woodrow \cite{W}) and finite directed graphs (Cherlin, \cite{Ch}). 

\subsection{Ramsey degrees} 

\label{subsection:Rdeg}

Having the Ramsey property is extremely restrictive for a class of finite structures. For that reason, weaker partition properties were introduced. One of the most common ones is obtained by imposing that colorings should only take a small number of colors on a large set, as opposed to being constant. This is captured by the following notion: for $k,l \in \N \smallsetminus \{ 0 \}$ and $L$-structures $\m{A}, \m{B}, \m{C}$, write $$\m{C} \arrows{(\m{B})}{\m{A}}{k,l}$$ when for any $c : \funct{\binom{\m{C}}{\m{A}}}{[k]}$ there is $\widetilde{\m{B}} \in \binom{\m{C}}{\m{B}}$ such $c$ takes at most $l$-many values on $\binom{\widetilde{\m{B}}}{\m{A}}$. Note that when $l = 1$, this is simply the partition property $\m{C} \arrows{(\m{B})}{\m{A}}{k}$ introduced previously. 

\begin{defn}

Let $\mathcal{K}$ be a class of $L$-structures. An element $\m{A} \in \mathcal{K}$ has a \emph{finite Ramsey degree in} $\K$ when there exists $l \in \N$ such that for any $\m{B} \in \mathcal{K}$, and any $k \in \N \smallsetminus \{ 0 \}$, there exists $\m{C} \in \mathcal{K} $ such that: $$\m{C} \arrows{(\m{B})}{\m{A}}{k,l}.$$

The least such number $l$ is denoted $\mathrm{t}_{\mathcal{K}}(\m{A})$ and is the \emph{Ramsey degree of} $\m{A}$ \emph{in} $\mathcal{K}$. 

\end{defn}

Equivalently, if $\mathcal K$ is Fra\"iss\'e and $\m{F}$ denotes its limit, $\m{A}$ has a finite Ramsey degree in $\mathcal{K}$ when there is $l \in \N$ such that for any $\m{B} \in \mathcal{K}$, and any $k \in \N \smallsetminus \{ 0 \}$, $$\m{F} \arrows{(\m{B})}{\m{A}}{k,l}.$$

The Ramsey degree is then equal to the least such number $l$. Note that it depends only on $\m A$ and $\K$. Finite Ramsey degrees can be seen in two different ways. They reflect the failure of the Ramsey property within a given class $\K$, but also reflect that arbitrary finite colorings can always be reasonably controlled.

As a concrete example, consider the class $\mathcal G$ of finite graphs. It is not Ramsey, but every $\m A \in \mathcal G$ has a finite Ramsey degree, which is equal to $$ t_{\mathcal G} (\m A) = |\m A|!/|\Aut(\m A)|.$$

\section{Compact flows}

\label{section:flows}

We now turn to topological groups and to dynamical properties of their actions. Let $G$ be a topological group. A \emph{$G$-flow} is a continuous action of $G$ on a topological space $X$. We will often use the notation $G \curvearrowright X$. The flow $G \curvearrowright X$ is \emph{compact} when the space $X$ is. It is \emph{minimal} when every $x \in X$ has dense orbit in $X$: \[ \forall x \in X \  \  \overline{G\cdot x} = X\] 

Finally, it is \emph{universal} when every compact minimal $G \curvearrowright Y$ is a factor of $G \curvearrowright X$, which means that there exists $\pi : X \longrightarrow Y$ continuous, onto, and so that $$\forall g \in G \quad \forall x \in X \quad  \pi (g \cdot x) = g \cdot \pi(x).$$ 

It turns out that when $G$ is Hausdorff, there is, up to isomorphism of $G$-flows, a unique $G$-flow that is both minimal and universal. This flow is called \emph{the universal minimal flow} of $G$, and is denoted $G \curvearrowright M(G)$. When the space $M(G)$ is reduced to a singleton, the group $G$ is said to be \emph{extremely amenable}. Equivalently, every compact $G$-flow $G\curvearrowright$ admits a fixed point, i.e. an element $x\in X$ so that $g\cdot x = x$ for every $g\in G$. We refer to \cite{KPT} or \cite{Pe} for a detailed account on those topics. Let us simply mention that, concerning extreme amenability, it took a long time before even proving that such groups exist, but that several non-locally compact transformation groups are now known to be extremely amenable (the most remarkable ones being probably the isometry groups of the separable infinite dimensional Hilbert space (Gromov-Milman, \cite{GM}), and of the Urysohn space (Pestov, \cite{Pe0})). As for universal minimal flows, prior to \cite{KPT}, only a few cases were known to be both metrizable and non-trivial, the most important examples being provided by the orientation-preserving homeomorphisms of the circle (Pestov, \cite{Pe1}), $S_{\infty}$ (Glasner-Weiss, \cite{GW1}), and the homeomorphism group of the Cantor space (Glasner-Weiss, \cite{GW2}). 

\section{The Kechris-Pestov-Todorcevic correspondence}

\label{section:KPT}

For an $L$-structure $\m A$, we denote by $\Aut (\m A)$ the corresponding automorphism group. Recall that when this group is trivial, we say that $\m A$ is rigid.

\begin{thm}[Kechris-Pestov-Todorcevic, \cite{KPT}, essentially Theorem 4.8]

\label{thm:EARP}

Let $\m F$ be a Fra\"{i}ss\'e structure, and let $G = \Aut(\m F)$. The following are equivalent: 

\begin{enumerate}
\item[i)] The group $G$ is extremely amenable. 

\item[ii)] The class $\Age (\m F)$ has the Ramsey property and consists of rigid elements.  
\end{enumerate}

\end{thm}

Because closed subgroups of $S_{\infty}$ are all of the form $\Aut(\m F)$, where $\m F$ is a Fra\"{i}ss\'e structure, the previous theorem actually completely characterizes those closed subgroups of $S_{\infty}$ that are extremely amenable. It also allows the description of many universal minimal flows via combinatorial methods. Indeed, when $\m F^{*} = (\m F, <^{*})$ is an order expansion of $\m F$, one can consider the space $\LO(\m F)$ of all linear orderings on $\m F$, seen as a subspace of $[2]^{\m F\times \m F}$. In this notation, the factor $[2]^{\m{F}\times \m F} = \{ 0, 1\}^{\m{F}\times \m F}$ is thought as the set of all binary relations on $\m F$. This latter space is compact, and $G$ continuously acts on it: if $S\in [2]^{\m{F}\times \m F}$ and $g\in G$, then $g\cdot S$ is defined by $$\forall x, y \in \m F \quad g\cdot S(x, y) \Leftrightarrow S(g^{-1}(x), g^{-1}(y)).$$ 

It can easily be seen that $\LO (\m F)$ and $X^{*} : = \overline{G \cdot <^{*}}$ are closed $G$-invariant subspaces.

\begin{thm}[Kechris-Pestov-Todorcevic, \cite{KPT}, Theorem 7.4]

\label{thm:OP}

Let $\m F$ be a Fra\" iss\'e structure in a language $L$, and $\m F^{*}$ a Fra\"iss\'e order expansion of $\m F$. The following are equivalent: 

\begin{enumerate}

\item[i)] The flow $G\curvearrowright X^{*}$ is minimal. 

\item[ii)] $\Age(\m F^{*})$ has the \emph{ordering property} relative to $\Age(\m F)$: for every $ \m A$ in $\Age (\m F)$, there exists $ \m B$ in $\Age(\m F)$ such that for every order expansion $\m A^{*}$ of $\m A$ in  $\Age(\m F^{*})$ and every order expansion $\m B^{*}$ of $\m B$ in $\Age (\m F^{*})$, $\m A^{*}$ embeds in $\m B^{*}$.

\end{enumerate}

\end{thm}

The following result, which builds on the two preceeding theorems, is then obtained: 

\begin{thm}[Kechris-Pestov-Todorcevic, \cite{KPT}, Theorem 10.8]

\label{thm:KPTUMF}

Let $\m F$ be a Fra\" iss\'e structure, and $\m F^{*}$ be a Fra\"iss\'e order expansion of $\m F$. The following are equivalent: 

\begin{enumerate}

\item[i)] The flow $G\curvearrowright X^{*}$ is the universal minimal flow of $G$. 

\item[ii)] The class $\Age(\m F^{*})$ has the Ramsey property as well as the ordering property relative to $\Age(\m F)$. 

\end{enumerate}

\end{thm}

A direct application of those results allowed to find a wealth of extremely amenable groups and of universal minimal flows. We will list many of the corresponding results later on, but let us mention at that point that some cases, which are very close to those described above, \emph{cannot} be captured directly by those theorems. Precisely, some Fra\"iss\'e classes do not have an order expansion with the Ramsey and the ordering property, but do so when the language is enriched with  additional symbols.

Therefore, we will not deal with order expansions in the language $L^{*}=L\cup \{<\}$ only (those will be later on referred to as \emph{pure} order expansions), but with \emph{precompact relational expansions}. For such expansions, we do not require $L^{*}=L\cup \{<\}$, but only $L^{*} = L\cup \{ R_{i}:i\in I\}$, where $I$ is countable, and every symbol $R_{i}$ is relational and not in $L$. An expansion $\m F^{*}$ of $\m F$ is then called \emph{precompact} when any $\m A \in \Age(\m F)$ only has finitely many expansions in $\Age(\m F^{*})$. Note that every $\m A \in \Age(\m F)$ has at least one expansion in $\Age(\m F ^{*})$: simply take a copy of $\m A$ in $\m F$, and consider the substructure of $\m F^{*}$ that it supports. The choice of the terminology is justified in \cite{NVT3}. For those expansions, the ordering property has a direct translation, which we call the \emph{expansion property}. 

\begin{defn}

Let $\K$ be a Fra\" iss\'e class in $L$ and let $\K^{*}$ be a relational expansion of $\mathcal K$. The class $\K^{*}$ has the \emph{expansion property} relative to $\K$ when for every  $ \m A \in \K$, there exists $ \m B \in \K$ such that $$\forall \m A^{*}, \m B^{*} \in \K^{*} \quad (\restrict{\m A^{*}}{L} = \m A \ \  \wedge \ \ \restrict{\m B^{*}}{L} = \m B) \Rightarrow \m A^{*} \leq \m B^{*}. $$
\end{defn}

Theorems \ref{thm:OP} and \ref{thm:KPTUMF} turn into the following versions: 

\begin{thm}

\label{thm:EP}

Let $\m F$ be a Fra\" iss\'e structure, and $\m F^{*}$ a precompact relational expansion of $\m F$ (not necessarily Fra\"iss\'e). The following are equivalent: 

\begin{enumerate}

\item[i)] The flow $G\curvearrowright X^{*}$ is minimal. 

\item[ii)] $\Age(\m F^{*})$ has the expansion property relative to $\Age(\m F)$. 

\end{enumerate}

\end{thm}

\begin{thm}

\label{thm:UMF}

Let $\m F$ be a Fra\" iss\'e structure, and $\m F^{*}$ be a Fra\"iss\'e precompact relational expansion of $\m F$. Assume that $\Age(\m F^{*})$ consists of rigid elements. The following are equivalent: 

\begin{enumerate}

\item[i)] The flow $G\curvearrowright X^{*}$ is the universal minimal flow of $G$. 

\item[ii)] The class $\Age(\m F^{*})$ has the Ramsey property as well as the expansion property relative to $\Age(\m F)$. 

\end{enumerate}

\end{thm}

%
%
%
%
%

One aspect should be emphasized here: the only reason for which the original paper \cite{KPT} was not written in the general setting we present here is that, at the time where it was developed, pure order expansions covered almost all known applications of the method to compute universal minimal flows (the cases that were left aside were computed easily with a bit of extra work). Arguably, they consequently constituted the right setting to establish a general correspondence. 

\section{Applications of the Kechris-Pestov-Todorcevic correspondence}

\label{section:applications}

In addition to its theoretic interest, the power of the KPT correspondence lies in its applications. Below is a list of Fra\"iss\'e classes for which it has been applied in order to compute the universal minimal flow of the corresponding automorphism groups. In all cases, the proof is combinatorial, and consists in finding a precompact expansion with the Ramsey and the expansion property (equivalently, an expansion where the number of expansions of each structure is equal to its Ramsey degree). 

\subsection{Graphs}

For simple, undirected, loopless graphs, the appropriate language is  $L = \{E\}$ with one binary relation symbol $E$. The symbol is then interpreted as a binary irreflexive and symmetric relation. The Fra\"iss\'e classes of finite graphs have been classified by Lachlan–Woodrow in \cite{LW}. In what follows, for each such class $\K$, we indicate a precompact expansion $\K^{*}$ which is Ramsey and has the expansion property:  

\begin{enumerate}

\item $\mathcal G$: all finite graphs. $\mathcal G^{*}$: pure order expansion consisting of all finite ordered graphs, i.e. add all the linear orderings (Abramson-Harrington \cite{AH} and Ne\v set\v ril-R\"odl \cite{NR1, NR2} independently). 

\item $\mathrm{Forb}(K_{n})$ ($n \geq 3$): $K_{n}$-free graphs, i.e. not containing the complete graph $K_{n}$ on $n$ vertices as a substructure. $\mathrm{Forb}(Kn)^{*}$: add all linear orderings. This result is due to Ne\v set\v ril-R\"odl \cite{NR1, NR2}, and was the first instance where the so-called \emph{partite construction} was used. This technique is still one of the most powerful tools in structural Ramsey theory. A simple account on it concerning triangle-free graphs can be found in \cite{GRS}. 

\item $\mathcal{EQ}$: disjoint unions of complete graphs, thought as the class of all finite equivalence relations. $\mathcal{EQ}$: add all \emph{convex} linear orderings (i.e. where all equivalence classes are intervals) (cf \cite{KPT}).  

\item $\mathcal{EQ}_{n}$ ($n\geq 1$): finite equivalence relations with at most n classes. $\mathcal{EQ}_{n}^{*}$: add unary relations symbols $P_{0},...,P_{n-1}$, interpreted as the parts of the equivalence relation, as well as all convex linear orderings which in a given structure $\m A$, order the parts as $P^{\m A}_{0}<^{\m A}...<^{\m A} P^{\m A} _{n-1}$ (essentially, Soki\'c \cite{So2}).  Note that if no requirement is put on the linear ordering, then the corresponding expansion is Ramsey (cf \cite{KPT}, p.158) but does not have the expansion property. 

\item $\overline{\mathcal{EQ}}_{n}$ ($n\geq 1$): finite equivalence relations, all of whose classes have at most $n$ elements. $\overline{\mathcal{EQ}}_{n}^{*}$: first, add all convex linear orderings (convex with respect to the equivalence relation). Then, add unary relations symbols $P_{0},...,P_{n-1}$, interpreted as disjoint transversals sets, and so that in each structure $\m A$, each $P^{\m A} _{i}$ has at most one element and satisfies $P^{\m A}_{0}<^{\m A}...<^{\m A} P^{\m A} _{n-1}$ (essentially, Soki\'c \cite{So2})

\item The complement of one of the classes $\K$ listed above, obtained by replacing all edges by non-edges and vice-versa. 
\end{enumerate}

\subsection{Hypergraphs}

For hypergraphs (also sometimes called set-systems), an appropriate language $L$ is made of countable many relational symbols with arity at least $2$. In a given structure $\m A$, each symbol $R$ of arity $n$ is interpreted as a relation $R^{\m A}$ so that if $R^{\m A}(a_{0},...,a_{n-1})$, then all $a_{i}$'s are distinct, and $R^{\m A}(a_{\sigma (0)},...,a_{\sigma (n-1)})$ for every permutation $\sigma$. Contrary to the case of graphs, no classification result of Fra\"iss\'e classes is available when $L$ is specified at this level of generality. Still, several results are known concerning hypergraphs. 

\subsubsection{General results}

Given $L$, the class $\mathcal H _{L}$ of \emph{all} finite hypergraphs in $L$ is a Fra\"iss\'e class, and one can take $\mathcal H _{L}$ to be the pure order expansion obtained by adding all linear orderings. This result extends the previous results about graphs, and has also been proved by independently by Abramson-Harrington and Ne\v set\v ril-R\"odl (same references as above). As in the case of graphs, substantial results have also been obtained by Ne\v set\v ril-R\"odl thanks to the partite construction when forbidden substructures are introduced. However, some technical requirements have to be satisfied by those forbidden substructures and we will not detail them here. 

\subsubsection{A particular case: boron tree structures}

Apart from those general results, some particular classes of hypergraphs were studied. Let us mention here the case of \emph{boron tree structures}, introduced by Cameron in \cite{Ca0}. Following \cite{Ca2}, let us say that a \emph{boron tree} is a finite (graph-theoretic) tree where all vertices have valency $1$ (``hydrogen atoms'') or $3$ (``boron atoms''). A boron tree $T$ gives raise to a \emph{boron tree structure} $\m T$, whose points are the hydrogen atoms and where $R^{\m T}(x,y,z,t)$ holds when the paths joining $x$ to $y$ and $z$ to $t$ do not intersect.  A precompact expansion with the Ramsey and the expansion property has been found by Jasi\'nski in \cite{J}, and was used to compute the corresponding universal minimal flow. Note that boron trees structures also provided one of the first examples where pure order expansions do not suffice in order to apply KPT theory. See also the paper \cite{Sol5} by Solecki for a different proof.

\subsection{Directed graphs}

Simple loopless directed graphs were classified by Cherlin in \cite{Ch}. The appropriate language here is made of one binary relation symbol $E$ interpreted in a structure $\m A$ as an irreflexive relation which also satisfies that at most one among $E^{\m A}(a,b)$ and $E^{\m A}(b,a)$ holds. We do not list here all the Fra\"iss\'e classes, but simply indicate several contributions that relate to the study of Ramsey expansions. 

\subsubsection{Posets}

The most important result appears in the paper \cite{PTW} by Paoli-Trotter -Walker, and deals with the class of all finite posets. It is shown there that the class of finite posets ordered by a linear extension has the Ramsey property. Note that this result is however attributed to Ne\v set\v ril-R\"odl (even though the corresponding paper was never published). The other Fra\"iss\'e classes of finite posets, which were classified in \cite{Sch}, are studied in detail by Soki\'c in \cite{So1, So2}. All the corresponding universal minimal flows can be found there. Note also that because of the correspondence between finite quasi-ordered sets and finite topological spaces, those results can be used in order to classify Ramsey classes of finite topological spaces (see \cite{Ne0}). 

\subsubsection{Local orders}

The paper \cite{LNS} details the Ramsey properties of another directed graph, the so-called dense local order. The corresponding results are presented in the language of precompact expansions in \cite{NVT3}. 

\

Finally, all the remaining cases of ultrahomogeneous directed graphs are treated in the recent work \cite{JLNW}.

\subsection{Metric spaces}

Even though metric spaces are very natural relational structures, it is only recently that they were studied from the point of view of structural Ramsey theory. An appropriate language here is $\{ d_{\alpha} : \alpha \in \Q\}$ where all the symbols are relational with arity $2$, and where, in a structure $\m A$, $d^{\m A} _{\alpha}(x,y)$ means that the distance is less than $\alpha$. The first significant result is due to Ne\v set\v ril in \cite{Ne}, where the class of all finite metric spaces is proved to be Ramsey when expanded with all linear orderings. This class not being countable, it is not Fra\"iss\'e, but it becomes so when the set of distances is restricted to some reasonable countable subset $S$ of the reals. When $S = \Q$ or $\N$, this allows to compute the universal minimal flow of the isometry groups of the rational Urysohn space, as well as of the usual Urysohn space. This was even the original motivation for Ne\v set\v ril's work, which came as a response to an early version of \cite{KPT} where it was asked whether a Ramsey theorem for metric spaces is hidden behind the extreme amenability of the isometry group of the Urysohn space (a result due to Pestov in \cite{Pe0}). Various other cases for $S$ are considered in \cite{NVT1}. This allows in particular to capture all Fra\"iss\'e classes of finite ultrametric spaces; for those, the relevant expansions are always with convex linear orderings. However, for most classes of metric spaces, the search for Ramsey expansions is largely open.

\subsection{Vector spaces}

Vector spaces over a finite field $F$ constitute a class for which the structural Ramsey properties were conjectured by Rota, and proved by Graham-Leeb-Rothschild in \cite{GLR1, GLR2} (see also \cite{Sp} for a shorter proof, based on Hales-Jewett theorem). The corresponding language consists of unary function symbols $M_{\alpha}$ ($\alpha \in F$) interpreted as scalar multiplications, and of a binary function symbol $+$ interpreted as addition. Those structures do not need to be expanded to become Ramsey, but they do (by antilexicographical linear orderings coming from an ordering of a basis) if one wants to apply the KPT correspondence.

\subsection{Lattices}

For lattices, the appropriate language is $L = \{ \land, \lor \}$, made of two binary function symbols, which are interpreted as the meet and the join operations respectively. The universal minimal flow was computed for the automorphism group of the random distributive lattice in \cite{KS}. The result uses in an essential way the Ramsey property for the class of Boolean lattices (distributive lattices equipped with a relative complementation). This latter result, also known as the dual Ramsey theorem, or Graham-Rothschild theorem (proved in \cite{GR}), is one of the basic tools of structural Ramsey theory. It can also be used to compute the universal minimal flow for the automorphism group of the coutable atomless Boolean algebra (for which the language is the previous language expanded by $\{-, 0, 1 \}$). Note also that all the corresponding groups are \emph{non-amenable} (for more on this topic, see \cite{KS}).

\subsection{Mixed structures}
Some other kinds of structures that also recently attracted some attention are those that are obtained by a superposition of several known simpler structures. For example, consider homogeneous permutations as defined by Cameron in \cite{Ca1}, and which correspond to structures equipped with two different linear orderings. The class of all homogeneous permutations is Ramsey (a result independently due to B\"ottcher-Foniok in \cite{BF} and Soki\'c in \cite{So0}). In fact, Fra\"iss\'e classes of homogeneous permutations have been classified in \cite{Ca1}, and their Ramsey properties are studied in \cite{BF}. Various other classes of structures that can be obtained by superposition were extensively studied by Soki\'c. Those include in particular structures equipped with several linear orderings and structures equipped with unary relations (see \cite{So3, So4, So5}). Finally, a general superposition theorem was obtained by Bodirsky in \cite{Bod} (see also \cite{So6} for a slightly different approach on this).

\section{Open questions and perspectives}

\label{section:open}

We close this article with a selection of open questions and perspectives related to the topics covered previously.  

\subsection{A general question concerning the existence of Ramsey expansions}

Among Fra\"iss\'e classes, a common point of view after the knowledge accumulated in the eighties is that Ramsey classes are quite exceptional objects. When analyzing how the most famous results of the field were obtained, it seems that two categories emerge. The first one corresponds to those ``natural'' classes where the Ramsey property holds: finite sets, finite Boolean algebras, finite vector spaces over a finite field. The second one corresponds to those classes where the Ramsey property fails but where this failure can be fixed by adjoining a linear ordering: finite graphs, finite $K_{n}$-free graphs, finite hypergraphs, finite partial orders, finite topological spaces, finite metric spaces... As for those classes where more than a linear ordering is necessary, besides the ones that appear in \cite{KPT} (finite equivalence relations with classes of size bounded by $n$, or equivalence relations with at most $n$ classes) or those that were found more recently and listed previously, not so many cases are known, but it would be extremely surprising that nobody encountered such instances before. Quite likely, the corresponding results were not considered as true structural Ramsey results, and were therefore overlooked. However, we saw that  precompact expansions seem to offer a reasonable general context, as they allow to compute the Ramsey degrees in the case of all Fra\"iss\'e classes of graphs, and directed graphs. In practice, it also appears that there is some sort of a standard scheme that can be applied in order to construct precompact Ramsey expansions whenever those exist. This motivates the following conjecture:  


\begin{conj}

\label{conj:precompactRamsey}
Let $\K$ be a Fra\"iss\'e class where there are only finitely many nonisomorphic structures in every cardinality (equivalently, $\K$ is the age of a countable ultrahomogeneous $\omega$-categorical structure). Then $\K$ admits a Ramsey precompact expansion. 
\end{conj}  

Of course, stating this problem as a conjecture and not as a question only reflects my own view, which is certainly heavily influenced by the large number of different combinatorial examples that are now available. At that point, it could very well be that all those examples are in fact very particular, for example because most of the corresponding languages are finite and binary. Note also even though the problem has now been circulating for some years, its first appearance in print is in the paper \cite{BPT} by Bodirsky, Pinsker and Tsankov (see Section 7). 

Contrary to the common opinion expressed at the top of the present paragraph, a positive answer would imply that after all, Ramsey classes are not so rare. Note also that there is another formulation of the conjecture. It is in terms of topological dynamics, and leaves open the possibility of a solution via techniques from dynamics and functional analysis. It also motivates the hypothesis made on $\K$, and shows that  the conjecture is false when no restriction is placed on $\K$.

\begin{thm}[Melleray-NVT-Tsankov \cite{MNT}]

\label{thm:metrizable}

Let $\m F$ be a Fra\"iss\'e structure, and let $G=\Aut(\m F)$. The following are equivalent: 

\begin{enumerate}

\item[i)] The structure $\m F$ admits a Fra\"iss\'e precompact expansion $\m F^{*}$ whose age has the Ramsey property and consists of rigid elements. 

\item[ii)] The flow $G\curvearrowright M(G)$ is metrizable and has a generic orbit.

\item[iii)] The group $G$ admits an extremely amenable closed subgroup $G^{*}$ such that the quotient $G/G^{*}$ is precompact. 

\end{enumerate} 

\end{thm}

%
%

Theorem \ref{thm:metrizable} is the reason for which Conjecture \ref{conj:precompactRamsey} is only made for Fra\"iss\'e classes where there are only finitely many non isomorphic structures in every cardinality. Indeed, there are many known closed subgroups of $S_{\infty}$ whose universal minimal flow is not metrizable (e.g. the countable discrete ones). Starting from those, the previous result produces some Fra\"iss\'e classes that do not have any precompact Ramsey expansion. For example, consider the structure $(\Z, d^{\Z}, <^{\Z})$ where $d^{\Z}$ and $<^{\Z}$ are the standard distance and ordering on $\Z$. Its automorphism group is $\Z$. Therefore, the corresponding age does not have any precompact Ramsey expansion (a fact which is actually easy to see directly).  

Theorem \ref{thm:metrizable} also allows to translate Conjecture \ref{conj:precompactRamsey} into purely dynamical terms. Call a closed subgroup of $S_{\infty}$ \emph{oligomorphic} when for every $n\in \N$, it induces only finitely many orbits on $\N^{n}$. Those groups are exactly the ones that appear as automorphism groups of Fra\"iss\'e structures whose age only has finitely many elements in every cardinality. Conjecture \ref{conj:precompactRamsey} then states that every closed oligomorphic subgroup of $S_{\infty}$ should have a metrizable universal minimal flow with a generic orbit. Using this terminology, the recent work of \cite {DGMR} by Dorais-Gubkin-McDonald-Rivera shows that in addition to all the groups coming from the aforementioned Fra\"iss\'e classes, the conjecture also holds for all the groups coming from $\omega$-categorical linear orders. At the moment, it is even possible that this should be true for a larger class of groups, called \emph{Roelcke precompact}. A topological group is such when it is precompact with respect to the greatest lower bound of the left and right uniformities. For a closed subgroup of $S_{\infty}$, being Roelcke precompact is equivalent to being an inverse limit of oligomorphic groups (see \cite{Ts}), and it turns out that so far, all known universal minimal flows coming from Roelcke precompact groups are metrizable with a generic orbit. 

It is natural to ask whether the assumption of the existence of the $G_{\delta}$ orbit is really necessary in item ii) of Theorem~\ref{thm:metrizable}. This question in fact appears in \cite{AKL} by Angel, Kechris, and Lyons, where it is asked whether for $G$ Polish, $M(G)$ necessarily has a $G_\delta$ orbit when it is metrizable. It turns out that when $G$ is a subgroup of $S_\infty$ (the most interesting case in the present survey), the answer is positive, and due to Zucker \cite{Z}. This is done via completely different techniques from \cite{MNT} by proving a stronger version of Theorem~\ref{thm:metrizable}, where a new item ii'), obtained from ii) by removing the hypothesis of the existence of a $G_\delta$ orbit, is added.

\subsection{Particular open problems in finite Ramsey theory}

Concerning particular open problems in finite Ramsey theory, we can mention finite metric spaces with distances in some set $S$, Euclidean metric spaces (this problem appears in \cite{KPT}), projective Fra\"iss\'e classes (those are developed in \cite{IS} and are connected to Fra\"iss\'e classes of finite Boolean algebras) and equidistributed Boolean algebras (this problem appears in \cite{KST}). In view of those problems, which very likely require the introduction of new techniques, some attention must be paid to the recent work \cite{Sol2}, \cite{Sol3} and \cite{Sol4} by Solecki, and \cite{T} by Todorcevic. The papers \cite{Sol2}, \cite{Sol3} reprove and generalize some of the classical results from structural Ramsey theory (see also \cite{Vl} for an infinite version of one of those results). On the other hand, both \cite{Sol4} and \cite{T} provide a unified approach that allows to derive the main basic results of Ramsey theory (classical Ramsey, Graham-Rothschild, Hales-Jewett) thanks to a unified abstract framework.

\subsection{Thomas conjecture}

Another question has to do with a conjecture of Thomas and was asked in the recent work of subsection and Pinsker \cite{BP}. Following \cite{BP}, let us say that a \emph{reduct} of a relational structure $\m A$ is a relational structure with the same domain as $\m A$ all of whose relations can be defined by a first-order formula in $\m A$. Thomas conjectured in \cite{Th} that every Fra\"iss\'e relational structure $\m F$ in a finite language only has finitely many reducts up to first-order interdefinability. Can anything be said if  $\Age(\m F)$ consists of rigid elements and has the Ramsey property? Note also that the recent progress concerning Thomas conjecture and the classification of reducts of classical ultrahomogeneous structures produced quite a number of new Fra\"iss\'e classes, for which the study of Ramsey properties may be quite accessible.

\subsection{Dynamics}

Besides the KPT correspondence, the paper \cite{KPT} has recently been related to quite a number of promising developments. Some take the correspondence to different contexts. It is the case for \cite{B} by Barto\v sov\'a on structures that are not necessarily countable. It is also the case for the paper \cite{MT} by Melleray-Tsankov who show how it can be transferred to the so-called metric Fra\"iss\'e structures. What is interesting here is that the equivalence between Ramsey property and extreme amenability of non-Archimedean Polish groups becomes an equivalence between an approximate version of the Ramsey property and extreme amenability of \emph{all} Polish groups. This equivalence actually captures some prior result obtained by Pestov in \cite{Pe0}, but because of the lack of technique to prove the approximate Ramsey property, it has not led to any practical result so far. Still, the parallel between classical and metric Fra\"iss\'e theory seems worth investigating. 

More generally, the combinatorial translation of dynamical facts performed in \cite{KPT} opens a variety of perspectives connected to combinatorics and dynamics. For example, the usual notion of amenability can actually be studied via two different approaches. The first one has to do with universal minimal flows, since a topological group is amenable if and only if its universal minimal flow admits an invariant Borel probability measure. As a direct consequence, $S_{\infty}$, the automorphism group of the countable random graph or the isometry group of the rational Urysohn space are amenable, but the automorphism groups of the countable atomless Boolean algebra or of the countable generic poset are not (see \cite{KS} by Kechris and Soki\'c). The second approach relative to amenability consists in expressing it directly in combinatorial terms using a ``convex'' version of the Ramsey property. This was done recently and by Moore in \cite{M} and by Tsankov (private communication). This approach did not lead to any concrete result so far, but most probably because nobody has really tried to develop techniques in direction of the convex Ramsey property. 

Amenability is also connected to another combinatorial condition called the Hrushovski poperty. A Fra\"iss\'e class $\mathcal K$ of finite structures satisfies the Hrushovski property when for every $\m A \in \mathcal K$, there exists $\m B \in \mathcal K$ containing $\m A$ so that every isomorphism between finite substructures of $\m A$ extends to an automorphism of $\m B$. It is proved by Kechris and Rosendal in \cite{KR} that the Hrushovski property translates nicely at the level of automorphism groups. Namely, it is equivalent to the fact that there is an increasing sequence of compact groups whose union is dense in $\Aut(\m F)$, where $\m F = \mathrm{Flim} (\mathcal K)$. Therefore, the Hrushovski property for $\mathcal K$ implies the amenability of $\Aut(\m F)$. It is also central in the study of other properties of Polish groups like the small index property, the automatic continuity property, and the existence of ample generics (see for example \cite{KR}, \cite{Sol1} or more recently \cite{Ke}). Nevertheless, there are still very natural classes of structures for which the Hrushovski property is not known to hold (e.g. the class of all finite tournaments) and those provide good, potentially difficult, combinatorial problems.  

Still in connection with amenability, the paper \cite{AKL} has recently pointed out an intriguing fact: every known case of non-Archimedean Polish group $G$ which is amenable and has a metrizable universal minimal flow turns out to be \emph{uniquely ergodic}, in the sense that every compact $G$-flow has a unique invariant Borel probablity measure (which is then necessarily ergodic). The question of knowing whether this always holds is open. Unique ergodicity also appears to be connected to new combinatorial phenomena, such as the uniqueness of so-called \emph{consistent random orderings} or a quantitative version of the expansion property. This last property can be described as follows: let $\m F$ be a Fra\"iss\'e structure and let $\m F^{*}$ be a Fra\"iss\'e precompact expansion of $\m F$. Say that $\Age(\m F^{*})$ has the \emph{quantitative expansion property} relative to $\Age(\m F)$ when there exists an isomorphism invariant map $\rho : \funct{\Age(\m F^{*})}{[0,1]}$ such that for every $\m A\in \Age(\m F)$ and every $\varepsilon >0$, there exists $\m B \in \Age(\m F)$ and a non-empty set of embeddings $E(\m A, \m B)$ of $\m A$ into $\m B$ such that for every expansion $\m A^{*}$ of $\m A$ and $\m B^{*}$ of $\m B$ in $\Age(\m F^{*})$, the proportion of embeddings in $E(\m A, \m B)$ that induce embeddings of $\m A^{*}$ in $\m B^{*}$ equals $\rho(\m A^{*})$ within $\varepsilon$.  Here is now how the quantitative expansion property is connected to unique ergodicity: if $\Age(\m F^{*})$ has the Ramsey property as well as the expansion property relative to $\Age(\m F)$, then the quantitative expansion property implies that $\Aut(\m F)$ is uniquely ergodic as soon as $\Aut(\m F)$ is amenable. This result, which is at the heart of \cite{AKL}, is used to show that several well-known automorphism groups are uniquely ergodic (e.g. of the random graph, of the ultrahomogeneous $K_{n}$-free graphs and of the rational Urysohn space).  

Last, it seems that a number of classical dynamical notions can be studied in the context of non-Archimedean Polish groups. 

Indeed, topological dynamics is full of many other natural classes of compact flows admitting universal minimal objects (e.g. equicontinuous flows, distal flows, almost periodic flows,...see \cite{dV}). Each of them is a potential candidate for an analog of \cite{KPT}, with potential new combinatorial and dynamical phenomena (For example, a strategy similar to \cite{KPT} can be applied in the context of proximal flows. The corresponding results appear in \cite{MNT}).

%
%

\bibliographystyle{amsalpha}
\bibliography{Bib14Dec}
\end{document}